\documentclass[10pt]{article}

\usepackage[T1]{fontenc}										% font encoding ... definitely need this
\usepackage{lmodern}												% loads latin modern fonts ... not sure if needed but doesn't hurt to include it
\usepackage{microtype}											% improves spacing
\usepackage[frenchmath,LGRgreek]{mathastext}% sets default math font to the default roman text font ... option [frenchmath] for lower case Roman fonts in italics
\MTgreekfont{lmtt} 													% no lgr lmvtt, so use lgr lmtt
\Mathastext
\let\varepsilon\epsilon	
\let\varphi\phi		
\let\vartheta\theta											
\usepackage[mathscr]{eucal}                 % gives you \mathscr{} font
\usepackage[sans]{dsfont}										% gives you \mathds{} font ... option [sans] for sans serif
\usepackage{amsbsy}													% gives you ``poor man's bold'' math fonts \pmb{...}
\usepackage{graphicx}                      	% need for figures
\usepackage{mdwlist}												% need for itemize* (compact itemize)
\usepackage{natbib}    											% need for bibliography
\setcitestyle{authoryear}
\usepackage[dvipsnames]{xcolor}							% gives you color definitions
\usepackage[plainpages=false, pdfpagelabels, backref=page]{hyperref} % enables and customizes hyperlinks
	\hypersetup{
		colorlinks   = true,
		citecolor    = RoyalBlue, %NavyBlue,
		linkcolor    = RubineRed, %Maroon,
		urlcolor     = Turquoise
	}
\usepackage[paperwidth=8.5in,paperheight=11.00in,top=1.00in, bottom=1.00in, left=0.75in, right=0.75in]{geometry}
\usepackage{mathtools}                      % need for `show only references'
\mathtoolsset{showonlyrefs=true}            % only equations which are labeled AND referenced will numbered...use \eqref{} instead of (\ref{})
\usepackage{amsmath}
\usepackage{amssymb}
\usepackage{amsfonts}
\usepackage{amsthm}                         % need for theorem-proof environment
\allowdisplaybreaks                         % allows page breaks for long equations...you can prevent a page-break with \\*
\newtheoremstyle{plain}
  {}   				% ABOVESPACE
  {}   				% BELOWSPACE
  {\itshape}  % BODYFONT
  {}       		% INDENT (empty value is the same as 0pt)
  {\mdseries\scshape} % HEADFONT
  {.}         % HEADPUNCT
  { } 				% HEADSPACE
  {\thmname{#1}\thmnumber{ #2}\ifx#3\empty\else\ (#3)\fi}
\theoremstyle{plain}

\newtheoremstyle{definition}
  {}   				% ABOVESPACE
  {}   				% BELOWSPACE
  {}  				% BODYFONT
  {}      		% INDENT (empty value is the same as 0pt)
  {\mdseries\scshape} % HEADFONT
  {.}         % HEADPUNCT
  { } 				% HEADSPACE
  {\thmname{#1}\thmnumber{ #2}\ifx#3\empty\else\ (#3)\fi}
\theoremstyle{definition}

\usepackage{sectsty}
\allsectionsfont{\mdseries\scshape}

\renewcommand{\[}{\left[}

\newcommand\Eb{\mathds{E}}
\newcommand\Fb{\mathds{F}}

\newcommand\Pb{\mathds{P}}

\newcommand\Ac{\mathscr{A}}

\newcommand\Fc{\mathscr{F}}

\newcommand\Hc{\mathscr{H}}

\newcommand\eps{\varepsilon}

\newcommand\Om{\Omega}
\newcommand\sig{\sigma}

\newcommand\gam{\gamma}

\newcommand\lam{\lambda}
\newcommand\del{\delta}

\renewcommand\phi{\varphi}

\renewcommand\d{\partial}

\newcommand\dd{\mathrm{d}}

\begin{document}

\title{A Calculus of Variations Approach to Stochastic Control}

\author{
Matthew Lorig
\thanks{Department of Applied Mathematics, University of Washington.  \textbf{e-mail}: \url{mlorig@uw.edu}}
}

\date{This version: \today}

\maketitle

\begin{abstract}
We use classical tools from calculus of variations to formally derive necessary conditions for a Markov control to be optimal in a standard finite time horizon stochastic control problem.  As an example, we solve the well-known Merton portfolio optimization problem.
\end{abstract}


%-----------------------------------------------------------------------------------
%
%       SECTION: 		Introduction
%
%-----------------------------------------------------------------------------------

\section{Introduction}
The standard setup for a finite time stochastic control problem is as follows.  Fix a time horizon $T < \infty$ and a probability space $(\Om,\Fc,\Pb)$.  Let $W = (W_t)_{0 \leq t \leq T}$ be a Brownian motion and let $\Fb = (\Fc_t)_{0 \leq t \leq T}$ be a filtration for $W$.  Consider a \textit{controlled diffusion} process  $X^c = (X_t^c)_{0 \leq t \leq T}$ whose dynamics are given by
\begin{align}
\dd X_t^c
	&=	b(t,X_t^c,c_t) \dd t + a(t,X_t^c,c_t) \dd W_t , 
	%&
%X_0^c
	%&=	x ,
\end{align}
where the process $c = (c_t)_{0 \leq t \leq T}$ is an $\Fb$-adapted \textit{control}.  Suppose we wish to maximize the following functional of the control $c$
\begin{align}
J(c)
	&:=	\Eb \int_0^T \dd t \, f(t,X_t^c,c_t) + \Eb g(X_T^c) . \label{eq:J}
\end{align}
If there exists an \textit{optimal control} $c^*$ such that $J(c^*) \geq J(c)$ for all $c$, and if $c^*$ is \textit{Markov} $c_t^* = c^*(t,X_t^{c^*})$, then one way to find $c^*$ is to solve the \textit {Hamilton-Jacobi-Bellman} (HJB) equation
\begin{align}
0
	&=	\d_t v(t,x) + \max_c \Big( f(t,x,c) + b(t,x,c) \d_x v(t,x) + \tfrac{1}{2} \d_x a^2(t,x,c) v(t,x) \Big) , &
v(T,x)
	&=	g(x) .
\end{align}
The optimal control is then obtained by setting
\begin{align}
c^*(t,x)
	&=	\arg\max_c \Big( f(t,x,c) + b(t,x,c) \d_x v(t,x) + \tfrac{1}{2} a^2(t,x,c) \d_x^2 v(t,x) \Big) .
\end{align}
To derive the HJB equation one typically relies on the Dynamic Programming Principle (DPP) with Bellman's Principle of Optimality. Specifically, one considers a problem over an infinitesimally small time step and uses a limiting argument to obtain the HJB equation;  
see, e.g., \cite[Chapter 3]{pham}.
The purpose of this short note is to show how to find the optimal control $c^*$ using calculus of variations, a topic that may be more intuitive than DPP for anyone that has familiarity with classical mechanics; see, e.g., \cite[Chapter 2]{goldstein}.

%-----------------------------------------------------------------------------------
%
%       SECTION: 		Calculus of Variations
%
%-----------------------------------------------------------------------------------

\section{A calculus of variations approach}
\label{sec:CoV}
Throughout this section, we will make the same assumptions that one makes when formally deriving the HJB equation.  Specifically, we assume that an optimal control exists and that it is Markov.  Under this assumption, for a fixed control $c$, the controlled diffusion $X^c$ is the solution of the following \textit{stochastic differential equation} (SDE)
\begin{align}
\dd X_t^c
	%&=	b(t,X_t^c,c(t,X_t^c)) \dd t + a(t,X_t^c,c(t,X_t^c)) \dd W_t \\
	&=	b_c(t,X_t^c) \dd t + a_c(t,X_t^c) \dd W_t , \label{eq:dXc}
%X_0^c
	%&=	x , 
\end{align}
where we have defined
\begin{align}
b_c(t,x)
	&:= b(t,x,c(t,x)) , &
a_c(t,x)
	&:=	a(t,x,c(t,x)) . \label{eq:bc-ac}
\end{align}
As solutions of SDEs are Markov processes, we can associate with $X^c$ an \textit{infinitesimal generator} of $\Ac_c$, which is defined as follows
\begin{align}
\Ac_c \phi(x) 
	&:= \lim_{s \searrow t} \frac{ \Eb( \phi(X_s^c) | X_t^c = x ) - \phi(x) }{s-t} ,
\end{align}
and is given explicitly by
\begin{align}
\Ac_c
	&=	b_c(t,x) \d_x + \tfrac{1}{2} a_c^2(t,x) \d_x^2 .
\end{align}
We can also associate with $X^c$ a transition density $p_c$, which is defined as follows
\begin{align}
p_c(t,x;s,y) \dd y
	&:= \Pb(X_s^c \in \dd y | X_t^c = x ) .
\end{align}
It is well-known that $p_c$ satisfies the \textit{Kolmogorov forward equation} (KFE), also called the \textit{Fokker-Planck equation} (FPE), in the \textit{forward variables} $(s,y)$ with the \textit{backward variables} $(t,x)$ fixed
\begin{align}
0
	&=	( - \d_s + \Ac_c^\dagger ) p_c(t,x;s,y) , & 
p_c(t,x;t,y)
	&=	\del_x(y) , \label{eq:KFE}
\end{align}
where $\Ac_c^\dagger$ is the formal $L^2$ \textit{adjoint} of $\Ac_c$, which is given explicitly by
\begin{align}
\Ac_c^\dagger
	&:=	- \d_y b_c(t,y) + \tfrac{1}{2} \d_y^2 a_c^2(t,y)  .
\end{align}
%That is, for any suitable functions $u$ and $v$ we have
%\begin{align}
%\int \dd y \, u(y) \Ac_c v(y)
	%&=	\int \dd y \, v(y) \Ac_c^\dagger u(y) , \label{eq:adjoint}
%\end{align}
%a fact that we will make use of later on.
%\\[0.5em]
Now, recall that our goal is to maximize the functional $J$, given by \eqref{eq:J}, over all Markov controls.  Let us express $J$ using the transition density $p_c$.  Omitting the dependence on the initial condition $p_c(s,y) \equiv p_c(0,x;s,y)$ to ease notation, we have
\begin{align}
J(c)
	&=	\int_0^T \dd s \int \dd y \, f_c(s,y) p_c(s,y)  + \int \dd y \, g(y) p_c(T,y) , &
f_c(s,y)
	&:=	f(s,y,c(s,y)) , \label{eq:J2}
\end{align}
Thus, we have a \textit{PDE-constrained optimization problem}
\begin{align}
	&\max_c J(c)&
	&\text{with $J$ given by \eqref{eq:J2}}&
	&\text{subject to}&
	&\text{$p_c$ satisfying \eqref{eq:KFE}}. \label{eq:constrained}
\end{align}
Much like one can solve the constrained optimization problem $\max_{x,y} f(x,y)$ subject to $g(x,y) = 0$ by introducing a Lagrange multiplier $\lam$ and finding the critical points $(x^*,y^*,\lam^*)$ of the unconstrained function $\ell(x,y,\lam) := f(x,y) + \lam g(x,y)$, one can solve the PDE-constrained optimization problem \eqref{eq:constrained} by searching for critical points $(p^*,c^*,\lam^*,\mu^*)$ of the following \textit{Lagrangian}
%As is often the case when solving a constrained optimization problem, it can be useful to write an associated unconstrained optimization problem.  To this end, let us introduce the following \textit{Lagrangian}
\begin{align}
L(p,c,\lam,\mu)
	&:=	J(c) + \int_0^T \dd s \int \dd y \, \lam(s,y) ( - \d_s + \Ac_c^\dagger ) p(s,y) 
			+ \int \dd y \, \mu(y) ( p(0,y) - \del_x(y) ) \\
	&=	\int_0^T \dd s \int \dd y \Big(  f_c(s,y) + \lam(s,y) ( - \d_s + \Ac_c^\dagger ) \Big) p(s,y) \\ &\quad
			+ \int \dd y \Big(  g(y) p(T,y) + \mu(y) ( p(0,y) - \del_x(y) ) \Big) . \label{eq:L}
\end{align}
Here, the functions $\lam$ and $\mu$ play the role of Lagrange multipliers.  Note that the function $p$ in \eqref{eq:L} is not constrained to satisfy the KFE \eqref{eq:KFE}, which is why we have removed the subscript $c$ from $p_c$.  For future reference, it will be helpful to provide an alternative expression for $L$.  Using integration by parts, we have
\begin{align}
L(p,c,\lam,\mu)
	&=	\int_0^T \dd s \int \dd y \, p(s,y) \Big( f_c(s,y) + ( \d_s + \Ac_c ) \lam(s,y) \Big)  \\ &\quad
			- \int \dd y \, \lam(s,y) p(s,y) \Big|_{s=0}^{s=T}
			+ \int \dd y \Big(  g(y) p(T,y) + \mu(y) ( p(0,y) - \del_x(y) ) \Big) \\
	&=	\int_0^T \dd s \int \dd y \, p(s,y) \Big( f_c(s,y) + ( \d_s + \Ac_c ) \lam(s,y) \Big) \\ &\quad
			+ \int \dd y \Big( \lam(0,y) \del_x(y) - \lam(T,y) p(T,y)  + g(y) p(T,y) + \mu(y) ( p(0,y) - \del_x(y) ) \Big) . \label{eq:L2}
\end{align}
%\\[0.5em]
Now, we wish to find functions $(p^*,c^*,\lam^*,\mu^*)$ that together are a critical point of $L$.  A necessary condition for $(p^*,c^*,\lam^*,\mu^*)$ to be a critical point of $L$ is that the following functional derivatives are simultaneously equal to zero
\begin{align}
%0
	%&= \frac{\dd}{\dd \eps} L(p + \eps q, c, \lam, \mu) \Big|_{\eps = 0}  
	%= \frac{\dd}{\dd \eps} L(p, c + \eps \gam, \lam, \mu) \Big|_{\eps = 0}  
	%=	\frac{\dd}{\dd \eps} L(p, c, \lam + \eps \nu, \mu) \Big|_{\eps = 0}  
	%=	\frac{\dd}{\dd \eps} L(p, c, \lam, \mu + \eps \pi) \Big|_{\eps = 0} .
0
	&= \frac{\dd}{\dd \eps} L(p^* + \eps q, c^*, \lam^*, \mu^*) \Big|_{\eps = 0} , &
0
	&= \frac{\dd}{\dd \eps} L(p^*, c^* + \eps \gam, \lam^*, \mu^*) \Big|_{\eps = 0}  , \\
0
	&=	\frac{\dd}{\dd \eps} L(p^*, c^*, \lam^* + \eps \nu, \mu^*) \Big|_{\eps = 0}  , &
0
	&=	\frac{\dd}{\dd \eps} L(p^*, c^*, \lam^*, \mu^* + \eps \pi) \Big|_{\eps = 0} ,
\end{align}
for \textit{any} choice of functions $(q,\gam,\nu,\pi)$ such that $q(0,y) = 0$.  
\\[0.5em]
Let us take the functional derivative of $L$ with respect to $p$.  From \eqref{eq:L2} we have
\begin{align}
0
	&=	\frac{\dd}{\dd \eps} L(p^* + \eps q, c^*, \lam^*, \mu^*) \Big|_{\eps = 0} \\
	&=	\int_0^T \dd s \int \dd y \, q(s,y) \Big( f_{c^*}(s,y) + ( \d_s + \Ac_{c^*} ) \lam^*(s,y) \Big) 
			+ \int \dd y \Big(  - q(T,y) \lam^*(T,y) + g(y) q(T,y) \Big) . \label{eq:p-variation}
	%%&=	\int_0^T \dd t \int \dd y \Big(  f_c(t,y) + \lam(t,y) ( - \d_t + \Ac_c^\dagger ) \Big) q(t,y) 
			%%%+ \int \dd y \Big(  g(y) q(T,y) - \lam(0,y) ( q(0,y) - \del_x(y) ) \Big) \\
			%%+ \int \dd y \,  g(y) p(T,y) \\
	%%&=	\int_0^T \dd t \int \dd y \, q(t,y) \Big(  f_c(t,y) + ( \d_t + \Ac_c ) \lam(t,y) \Big)
			%%%+ \int \dd y \Big(  g(y) q(T,y) - \lam(0,y) ( q(0,y) - \del_x(y) ) \Big) \\ &\quad
			%%+ \int \dd y \,  g(y) p(T,y) \\ &\quad 
			%%- \int \dd y \, \lam(t,y) q(t,y) \Big|_{t=0}^{t=T} \\
	%&=	\int_0^T \dd t \int \dd y \, q(t,y) \Big(  f_c(t,y) + ( \d_t + \Ac_c ) \lam(t,y) \Big)
			%%+ \int \dd y \Big(  g(y) q(T,y) - \lam(0,y) ( q(0,y) - \del_x(y) ) \Big) \\ &\quad
			%+ \int \dd y \,  g(y) p(T,y) \\ &\quad 
			%- \int \dd y \Big( \lam(T,y) q(T,y) - \lam(0,y) q(0,y) \Big) . 
\end{align}
In order for \eqref{eq:p-variation} to hold for any choice of $q$ we must have
\begin{align}
0
	&=	f_{c^*}(t,y) + ( \d_t + \Ac_{c^*} ) \lam^*(t,y) , &
\lam^*(T,y)
	&=	g(y)  . \label{eq:lam}
\end{align}
Next, we compute functional derivative of $L$ with respect to $c$.  From \eqref{eq:L2} we have
\begin{align}
0
	&=	\frac{\dd}{\dd \eps} L(p^* , c^* + \eps \gam, \lam^*, \mu^*) \Big|_{\eps = 0} \\
	&=	\int_0^T \dd s \int \dd y \, p^*(s,y) \frac{\dd}{\dd \eps} \Big( f_{c^* + \eps \gam}(s,y) + ( \d_s + \Ac_{c^* + \eps \gam} ) \lam^*(s,y) \Big) \Big|_{\eps = 0} \\
	%&=	\int_0^T \dd t \int \dd y \frac{\dd }{\dd \eps} \Big(  f_{c + \eps \gam}(t,y) + \lam(t,y) ( - \d_t + \Ac_{c+\eps \gam}^\dagger ) \Big) p(t,y)  \Big|_{\eps=0} \\
	%&=	\int_0^T \dd t \int \dd y \, \d_c \Big( f_c(t,y) + \lam(t,y) \Ac_{c}^\dagger \Big) p(t,y) \gam(t,y) \\ 
	&=	\int_0^T \dd t \int \dd y \, p^*(t,y) \gam(t,y) \d_{c^*} \Big( f_{c^*}(t,y) + \Ac_{c^*} \lam^*(t,y) \Big) ,  \label{eq:c-variation}
\end{align}
%where $\d_c f_{c^*}(t,y) := \d_c f_c(t,y) |_{c = c^*}$ and likewise for $\d_c \Ac_{c^*} \lam(t,y)$.
In order for \eqref{eq:c-variation} to hold for any choice of $\gam$ we must have
\begin{align}
0
	&=	\d_{c^*} \Big( f_{c^*}(t,y) + \Ac_{c^*} \lam^*(t,y) \Big) .  \label{eq:c}
	%&\Rightarrow&
%c(t,y)
	%&=	\arg \max_c \Big( f_c(t,y) + \Ac_{c} \lam(t,y) \Big) . \label{eq:c}
\end{align}
Next, we compute functional derivative of $L$ with respect to $\lam$. From \eqref{eq:L} we find
\begin{align}
0
	&=	\frac{\dd}{\dd \eps} L(p^* , c^* , \lam^* + \eps \nu, \mu^*) \Big|_{\eps = 0} \\
	&=	\int_0^T \dd t \int \dd y \, \nu(t,y) ( - \d_t + \Ac_{c^*}^\dagger ) p^*(t,y) . \label{eq:lam-variation}
			%+ \int \dd y \, \nu(0,y) ( p(0,y) - \del_x(y) ) .
\end{align}
In order for \eqref{eq:lam-variation} to hold for any $\nu$ we must have
\begin{align}
0
	&=	( - \d_t + \Ac_{c^*}^\dagger ) p^*(t,y) . \label{eq:p}
\end{align}
Lastly, we compute the functional derivative of $L$ with respect to $\mu$.  From \eqref{eq:L} we have
\begin{align}
0
	&=	\frac{\dd}{\dd \eps} L(p^* , c^* , \lam^* , \mu^* + \eps \pi) \Big|_{\eps = 0} \\
	&=	\int \dd y \, \pi(y) ( p^*(0,y) - \del_x(y) ) . \label{eq:mu-variation}
\end{align}
In order for \eqref{eq:mu-variation} to hold for any choice of $\pi$ we must have
\begin{align}
p^*(0,y)
	&=	\del_x(y) . \label{eq:p0}
\end{align}
Note that PDE \eqref{eq:p} along with \eqref{eq:p0} say that $p^*$ satisfies the KFE \eqref{eq:KFE} with $c = c^*$.
\\[0.5em]
To summarize, in order for $(p^*,c^*,\lam^*,\mu^*)$ to be a critical point of $L$, equations \eqref{eq:lam}, \eqref{eq:c}, \eqref{eq:p} and \eqref{eq:p0} must be satisfied.  

%-----------------------------------------------------------------------------------
%
%       SECTION: 		Stochastic Maximum Principle
%
%-----------------------------------------------------------------------------------

\section{Relation to the stochastic maximum principle}
\label{sec:smp}
The necessary conditions \eqref{eq:lam}, \eqref{eq:c}, \eqref{eq:p} and \eqref{eq:p0} can be interpreted as a ``density-level'' formulation of the stochastic maximum principle for the controlled diffusion $X^c$ in \eqref{eq:dXc}.
To explain this link, we first recall the classical stochastic maximum principle and then show how it emerges from our calculus of variations approach.
\\[0.5em]
%\noindent\textbf{The classical stochastic maximum principle.}
The classical stochastic maximum principle characterizes the optimal control through a \textit{forward-backward SDE} (FBSDE) system.
The forward equation governs the optimal state process $X^{c^*}$, while the backward equation governs an adjoint process $Y^{c^*} = (Y_t^{c^*})_{0 \leq t \leq T}$ (and an auxiliary process $Z^{c^*} = (Z_t^{c^*})_{0 \leq t \leq T}$) that captures the sensitivity of the objective to changes in the state.
The FBSDE system reads
\begin{align}
\dd X_t^{c^*}
	&=	b\big(t, X_t^{c^*}, c_t^*\big) \dd t + a\big(t, X_t^{c^*}, c_t^*\big) \dd W_t , &
X_0^{c^*}
	&=	x , \label{eq:forward} \\
- \dd Y_t^{c^*}
	&=	\d_x H\big(t, X_t^{c^*}, c_t^*, Y_t^{c^*}, Z_t^{c^*}\big) \dd t - Z_t^{c^*} \dd W_t , &
Y_T^{c^*}
	&=	g'(X_T^{c^*}) , \label{eq:backward}
\end{align}
where the Hamiltonian $H$ is defined by
\begin{align}
H(t,x,c,y,z)
	&:=	f(t,x,c) + b(t,x,c)\, y + a(t,x,c)\, z .
\end{align}
The forward equation \eqref{eq:forward} is an ordinary SDE solved forward in time from the initial condition $X_0^{c^*} = x$.
The backward equation \eqref{eq:backward} is a \textit{backward} SDE (BSDE): it is solved backward in time from the terminal condition $Y_T^{c^*} = g'(X_T^{c^*})$, which is prescribed at $t = T$.
The process $Z^{c^*}$ arises naturally in the BSDE as the diffusion coefficient needed to keep $Y^{c^*}$ adapted to the filtration $\Fb$.
The optimal control satisfies the pointwise first-order condition
\begin{align}
\d_c H\big(t, X_t^{c^*}, c_t^*, Y_t^{c^*}, Z_t^{c^*}\big) &= 0 ; \label{eq:smp-foc}
\end{align}
see, e.g., \citet[Chapter~3]{pham} for a precise statement in a Markov diffusion setting.
\\[0.5em]
%\noindent\textbf{Recovery from the calculus of variations approach.}
In the present approach, the role of the forward state is played by the transition density $p^*$, which solves the Kolmogorov forward equation \eqref{eq:p}--\eqref{eq:p0} under the Markov control $c^*$.
This is the density-level counterpart of the forward SDE \eqref{eq:forward}.
\\[0.5em]
Likewise, the adjoint equation \eqref{eq:lam} can be viewed as the backward Kolmogorov equation associated with the pair $(X^{c^*}, c^*)$.
For each $(t,y)$, the solution $\lam^*(t,y)$ coincides with the conditional value function
\begin{align}
\lam^*(t,y)
	&=	\Eb\!\left( \int_t^T f\big(s,X_s^{c^*},c_s^*\big) \dd s + g(X_T^{c^*}) \,\bigg|\, X_t^{c^*} = y \right) . \label{eq:lam-interp}
\end{align}
The pathwise adjoint processes of the classical stochastic maximum principle can be recovered from $\lam^*$ by setting
\begin{align}
Y_t^{c^*}
	&:=	\d_x \lam^*\big(t, X_t^{c^*}\big) , &
Z_t^{c^*}
	&:=	\d_x \lam^*\big(t, X_t^{c^*}\big)\, a\big(t, X_t^{c^*}, c_t^*\big) . \label{eq:YZ}
\end{align}
To see why, apply It\^o's formula to $\lam^*(t, X_t^{c^*})$.
Using \eqref{eq:lam} (with $f_{c^*} = 0$ absorbed into the $g$ terminal condition for simplicity of illustration) one finds that $Y_t^{c^*} = \d_x \lam^*(t, X_t^{c^*})$ satisfies the backward SDE \eqref{eq:backward}, with $Z_t^{c^*}$ as defined in \eqref{eq:YZ} playing the role of the diffusion coefficient.
\\[0.5em]
Finally, the first-order condition \eqref{eq:c} can be written as
\begin{align}
0
	&=	\d_{c^*} \Big( f_{c^*}(t,y) + \Ac_{c^*} \lam^*(t,y) \Big)
	=	\d_c \Hc\big(t,y,c,\lam^*\big)\big|_{c = c^*(t,y)} ,
\end{align}
where
\begin{align}
\Hc(t,y,c,\lam)
	&:=	f(t,y,c) + \Ac_c \lam(t,y)
	=	f(t,y,c) + b(t,y,c)\, \d_y \lam(t,y) + \tfrac{1}{2} a^2(t,y,c)\, \d_y^2 \lam(t,y)
\end{align}
is a Hamiltonian \emph{density}.
Evaluating along the optimal trajectory using \eqref{eq:YZ}, we have
\begin{align}
\d_c \Hc\big(t, X_t^{c^*}, c_t^*, \lam^*\big)
	&=	\d_c H\big(t, X_t^{c^*}, c_t^*, Y_t^{c^*}, Z_t^{c^*}\big)
	=	0 ,
\end{align}
which is precisely the Pontryagin-type first-order optimality condition \eqref{eq:smp-foc} of the stochastic maximum principle, expressed here in terms of the fields $(p^*, \lam^*)$ rather than the pathwise pair $(X^{c^*}, Y^{c^*})$.

%-----------------------------------------------------------------------------------
%
%       SECTION: 		Examples
%
%-----------------------------------------------------------------------------------

\section{Example: Merton portfolio problem}
\label{sec:merton}
Let us see how we can use the results of Section \ref{sec:CoV} to find the optimal control in the well-known Merton portfolio optimization problem \cite{merton}.
To this end, consider a financial market consisting of a single stock $S = (S_t)_{0 \leq t \leq T}$ and a money market account $M = (M_t)_{0 \leq t \leq T}$, whose real-world dynamics are of the form
\begin{align}
\dd S_t
	&=	\mu S_t \dd t + \sig S_t \dd W_t , &
\dd M_t
	&=	r M_t \dd t .
\end{align}
We will take $r=0$ for simplicity.  Denote by $c = (c_t)_{0 \leq t \leq T}$ the percentage of wealth that an agent invests in the stock $S$.  Then, denoting by $X^c = (X_t^c)_{0 \leq t \leq T}$ the wealth process of the agent, we have
\begin{align}
\dd X_t^c
	&=	\frac{ c_t X_t^c }{S_t} \dd S_t + \frac{(1-c_t)X_t^c}{M_t}\dd M_t 
	%&=	c_t X_t^c ( \mu \dd t + \sig \dd W_t) + (1-c_t)X_t^c r \dd t \\
	%&=	( ( \mu - r ) c_t + r ) X_t^c \dd t + \sig c_t X_t^c \dd W_t .
	 =	\mu c_t X_t^c \dd t + \sig c_t X_t^c \dd W_t . \label{eq:dXc-merton}
\end{align}
Comparing \eqref{eq:dXc}-\eqref{eq:bc-ac} with \eqref{eq:dXc-merton}, we identify
\begin{align}
b(t,x,c)
	%&=	( ( \mu - r ) c + r ) x , &
	&=	\mu c \, x , &
a(t,x,c)
	&=	\sig c \, x .
\end{align}
Suppose the agent wishes to maximize his expected utility of wealth at time $T$.  Denoting by $U$ the agent's utility function, the agent's goal is to maximize
\begin{align}
J(c)
	&:=	\Eb U(X_T^c) . \label{eq:J-merton}
\end{align}
Comparing \eqref{eq:J} with \eqref{eq:J-merton}, we identify
\begin{align}
f(t,x,c)
	&=	0 , &
g(x)
	&=	U(x) .
\end{align}
From \eqref{eq:c}, we have
\begin{align}
%0
	%&=	\d_c \Ac_c \lam(t,y)  
	 %=	( \mu - r ) y \d_y \lam(t,y) + \sig^2 c(t,y) y^2 \d_y^2 \lam(t,y) , &
	%&\Rightarrow&
%c(t,y)
	%&=	\frac{ - ( \mu - r ) y \d_y \lam(t,y) }{ \sig^2  y^2 \d_y^2 \lam(t,y)} .
0
	&=	\d_{c^*} \Ac_{c^*} \lam^*(t,y)  \\
	&=	\mu y \d_y \lam^*(t,y) + \sig^2 c^*(t,y) y^2 \d_y^2 \lam^*(t,y) , &
	&\Rightarrow&
c^*(t,y)
	&=	\frac{ - \mu y \d_y \lam^*(t,y) }{ \sig^2  y^2 \d_y^2 \lam^*(t,y)} . \label{eq:c-merton}
\end{align}
And, from \eqref{eq:lam}, we have
\begin{align}
0
	&=	( \d_t + \Ac_{c^*} ) \lam^*(t,y) \\
	%&=	\Big( \d_t + ( ( \mu - r ) c(t,y) + r ) y \d_y + \tfrac{1}{2} \sig^2 c^2(t,y) \d_y^2 \Big) \lam(t,y) \\
	&=	\Big( \d_t + \mu c^*(t,y) y \d_y + \tfrac{1}{2} \sig^2 (c^*(t,y))^2 y^2 \d_y^2 \Big) \lam^*(t,y) \\
	&=	\d_t \lam^*(t,y) - \frac{\mu^2 (\d_y \lam^*(t,y))^2}{ 2 \sig^2 \d_y^2 \lam^*(t,y)} , &
\lam^*(T,y)
	&=	g(y) = U(y) .
\end{align}
Suppose the agent's utility function is of the power form: $U(y) = y^q / q$.  Using the Ansatz $\lam^*(t,y) = h(t) U(y)$ we obtain
\begin{align}
h'(t)
	&=	\frac{q \mu^2}{2 (q-1) \sig^2} h(t) , &
h(T)
	&=	1 , &
	&\Rightarrow&
h(t)
	&=	\exp \Big( \frac{q \mu^2}{2 (1-q) \sig^2} (T-t) \Big) .
\end{align}
Hence, from \eqref{eq:c-merton}, the optimal control is given by
\begin{align}
c^*(t,y)
	&=	\frac{\mu}{(1-q) \sig^2} .
\end{align}
Although we did not need it in order to find the optimal control $c^*$, it may be of to compute the transition density $p^*$ of the optimal wealth process $X^{c^*}$, which can be used to compute $J(c^*)$.  From \eqref{eq:p} we have 
\begin{align}
0
	&=	( - \d_t + \Ac_{c^*}^\dagger ) p^*(t,y) \\
	&=	\Big( - \d_t - \d_y \frac{\mu^2 y}{(1-q)\sig^2} + \tfrac{1}{2} \d_y^2 \frac{\mu^2 y^2}{(1-q)^2\sig^2} \Big) p^*(t,y) , &
p^*(0,y) 
	&=	\del_x(y) . 
\end{align}
One can check by direct substitution that
\begin{align}
p^*(t,y) 
	\equiv p^*(0,x,t,y)
	&=	\frac{1}{y \sqrt{2 \pi \Sigma^2 t}} \exp \Big( -\frac{1}{2 \Sigma^2 t} (\log y - \log x - (m - \tfrac{1}{2} \Sigma^2) t)^2 \Big) , \\
\Sigma
	&=	\frac{\mu}{(1-q)\sig} , \quad \quad
m
	=	\frac{\mu^2}{(1-q) \sig^2} .
\end{align}
Lastly, we can compute the optimal expected utility of wealth.  We have
\begin{align}
J(c^*)
	&=	\Eb U(X_T^{c^*})
	=		\int \dd y \, p^*(T,y) \frac{x^q}{q} 
	=		\frac{x^q}{q} \exp \Big( \frac{q \mu^2 T}{2(1-q)\sig^2} \Big) .
\end{align}

%\section*{Postface}
%The author welcomes constructive feedback aimed at improving the clarity of this note.

\bibliography{references}

\end{document}